\documentclass[leqno,12pt]{article} 
\usepackage{amsmath, amssymb}
\usepackage{amsthm} 
\theoremstyle{plain} 
\newtheorem{theorem}{\indent\sc Theorem}[section] 
\newtheorem{lemmama}[theorem]{\indent\sc Lemma}
\newtheorem{corollaryollary}[theorem]{\indent\sc Corollary}

\theoremstyle{definition} 
\newtheorem{definition}[theorem]{\indent\sc Definition}
\newtheorem{remark}[theorem]{\indent\sc Remark}

\makeatletter
\def\address#1#2{\begingroup
\noindent\parbox[t]{7.8cm}{%
\small{\scshape\ignorespaces#1}\par\vskip1ex
\noindent\small{\itshape E-mail address}%
\/: #2\par\vskip4ex}\hfill%
\endgroup}%
\makeatother
\title{\uppercase{Projective normality of 
  toric 3-folds  \\
 with non-big adjoint hyperplane sections, II}} 
\author{
\textsc{Shoetsu Ogata$^{*}$ } 
}
\date{} 
\begin{document}

\maketitle

\footnote{ 
2000 \textit{Mathematics Subject Classification}.
Primary 14M25; Secondary 52B20.
}
\footnote{ 
\textit{Key words and phrases}. 
Toric varieties, projectively normal, Fano variety.
}
\footnote{ 
$^{*}$Partly supported by the Grant-in-Aid for Scientific Research (C),
Japan Society for the Promotion of Science. 
}

\begin{abstract}
Let $(X, A)$ be  a nonsingular polarized toric 3-fold.
We show that if the adjoint  bundle  of $A$  has no global sections, then 
all ample line bundles  on $X$
 are  normally generated.  Even if the adjoint bundle of $A$ is effective, if it is not big, then
 it is also shown the normal generation. 
Especially, we show that all ample line bundles on  a nonsingular toric Fano 3-fold are
normally generated.
\end{abstract}

\section*{Introduction} 

An ample line bundle $L$ on a projective variety is called {\it  normally generated}
by Mumford \cite{Mf},
 if the multiplication map $\Gamma(L)^{\otimes i} \to
\Gamma(L^{\otimes i})$ is surjective for all $i\ge1$.  
If  $L$ is normally generated, then 
 it is very ample.   Furthermore, if the variety $X$ is normal,
then a normally generated ample line bundle $L$ defines the embedding
$\Phi_L: X \to \mathbb{P}(\Gamma(L))$ of $X$ as a {\it projectively normal}
variety, i.e., the homogeneous coordinate ring is a normal ring.  

For an ample line bundle $L$ on a toric variety of dimension $n$ ($>1$), we see that
the twisted bundle $L^{\otimes k}$ is normally generated for $k\ge n-1$ (see \cite{K} for $n=2$, in general case \cite{EW}  and \cite{N}).
For more detail, see the introduction of \cite{Og2}.

Ogata \cite{Og2}  proved the following theorem.
 
 \begin{theorem}[Ogata]\label{int:tm0}
 Let $L$ be an ample line bundle on a nonsingular projective toric variety $X$ of dimension three.  
 \begin{itemize}
 \item[{\rm (1)}] If $H^0(X, L\otimes\mathcal{O}_X(K_X))=0$,
 then $L$ is normally generated.
 
 \item[{\rm (2)}] If $H^0(X, L\otimes\mathcal{O}_X(K_X))\not=0$ and if $L\otimes\mathcal{O}_X(K_X)$ 
 is not big, then  $L$ is normally generated.
 \end{itemize}
 \end{theorem}
 
 In this article, we will prove the following theorems.
 
  \begin{theorem}\label{int:t1}
Let $(X,A)$ be a polarized nonsingular  toric variety of dimension three.  If $H^0(X, A\otimes\mathcal{O}_X(K_X))=0$,
 then any ample line bundle $L$ on $X$ is normally generated.
\end{theorem}

  \begin{theorem}\label{int:t2}
Let $(X,A)$ be a polarized nonsingular  toric variety of dimension three.  If $H^0(X, A\otimes\mathcal{O}_X(K_X))\not=0$ and if $A\otimes\mathcal{O}_X(K_X)$ is not big,
 then any ample line bundle $L$ on $X$ is normally generated.
\end{theorem}

\begin{corollaryollary}
Any ample line bundle on a nonsingular toric Fano variety of dimension three is
normally generated.
\end{corollaryollary}

In order to prove theorems, we use the following theorem of Ogata \cite{Og1}.
\begin{theorem}\label{keytm}
Let $X$ be a nonsingular projective toric variety of dimension three with a non-trivial
morphism onto the projective line.
Then all ample line bundles on $X$ are normally generated.
\end{theorem}

Since the statements in the theorems can be interpreted into equalities on lattice points of
 integral convex polytopes,
 for a proof of Theorems~\ref{int:t1} and \ref{int:t2}
we need to investigate  properties of convex polytopes of dimension three.

The structure of this paper is as follows:

In Section 1  we recall basic results about toric varieties and line bundles on them.

In Section 2 we give a proof of Theorem~\ref{int:t1}
by using the Ogata's classification\cite[Proposition 2.3]{Og2}
 of nonsingular integral convex polytopes
of dimension three without interior lattice points.

In Section 3 we first give a coarse characterization of lattice polytopes with small
internal polytopes.   By using this characterization we prove Theorem~\ref{int:t2}.


\section{Projective toric varieties}\label{sect1}

In this section, we recall the facts on toric varieties which we need in this paper following Oda's
book \cite{Od}, or Fulton's book \cite{Fu}.
For simplicity, we consider toric varities are defined over the complex number field.

Let $N$ be a free $\mathbb{Z}$-module of rank $n$, $M$ its dual and 
$\langle, \rangle: M\times N \to\mathbb{Z}$ the canonical pairing.  
By the scalar extension to the field $\mathbb{R}$ of real
numbers, we have real vector spaces $N_{\mathbb{R}}:=N\otimes_{\mathbb{Z}}\mathbb{R}$
and $M_{\mathbb{R}}:= M\otimes_{\mathbb{Z}}\mathbb{R}$.  We denote also by $\langle, \rangle$
 the pairing of $M_{\mathbb{R}}$
and $N_{\mathbb{R}}$ defined by the scalar extension.
Let $T_N:= N \otimes_{\mathbb{Z}}\mathbb{C}^* \cong (\mathbb{C}^*)^n$ be the algebraic 
torus over the field $\mathbb{C}$ of complex numbers, where $\mathbb{C}^*$ is the multiplicative
group of $\mathbb{C}$.  Then the character group 
$\mbox{Hom}_{\scriptstyle{gr}}(T_N, \mathbb{C}^*)$ of $T_N$ is identified with $M$
 and $T_N=\mbox{Spec}\ \mathbb{C}[M]$.
For $m\in M$ we denote $\bold{e}(m)$ as the character of $T_N$.  Let $\Delta$ be a finite complete
fan in $N$ and  $X=T_N \mbox{emb}(\Delta)$ a complete toric variety
 of dimension $n$ (see  \cite[Section 1.2]{Od}, or  \cite[Section 1.4]{Fu}).
  We note that a toric variety defined by a fan is always normal.

Let $L$ be an ample  line bundle on $X$.  Then
we have an integral convex polytope $P$ in $M_{\mathbb{R}}$ with
\begin{equation}\label{sect1:eq2}
H^0(X, L) \cong \bigoplus_{m\in P\cap M}\mathbb{C}\bold e(m),
\end{equation}
where $\bold e(m)$ are considered as rational functions on $X$ because
they are functions on an open dense subset $T_N$ of $X$ (see 
  \cite[Section 2.2]{Od}, or \cite[Section 3.5 ]{Fu}).  
Here an integral convex polytope $P$ in $M_{\mathbb{R}}$ is the convex
hull $\mbox{Conv}\{u_1, u_2, \dots, u_s\}$ in $M_{\mathbb{R}}$ of
a finite subset $\{u_1, u_2, \dots, u_s\} \subset M$.
We note that $\dim_{\mathbb{R}} P=\dim X$.
The $l$-th power  $L^{\otimes l}$ corresponds to 
the convex polytope
$lP :=\{lx\in M_{\mathbb{R}}; x\in P\}$.

\begin{definition}
An integral convex polytope $P$ in $M_{\mathbb{R}}$ of dimension $n$ is
called {\it nonsingular} if for each vertex $u$ of $P$ the cone $\mathbb{R}_{\ge0}(P-u)
:=\{\lambda(x-u_i)\in \mathbb{R}^n;
x\in P\ \mbox{and} \ \lambda \ge0\}$
is nonsingular, that is, there exists a $\mathbb{Z}$-basis $\{m_1, \dots, m_n\}$ of $M$
such that
$$
\mathbb{R}_{\ge0}(P-u) = \mathbb{R}_{\ge0}m_1 + \dots + \mathbb{R}_{\ge0}m_n.
$$
\end{definition}

\begin{remark}   Set $V=\{v_1, \dots, v_r\}$  the set of all vertices of $P\subset M_{\mathbb{R}}$
with $\dim_{\mathbb{R}}P=\mbox{rank}\ M=n$.
Let $\sigma_i\subset N_{\mathbb{R}}$ be the dual cone to $\mathbb{R}_{\ge0}(P-v_i)$ with respect to
the pairing $\langle, \rangle$.  Then the set of all faces of $\sigma_1, \dots, \sigma_r$
coincides with the fan $\Delta$ defining the toric variety $X$,
and we can define the ample line bundle $L$ on $X$ satisfying (\ref{sect1:eq2}).
See \cite[Section 2.4]{Od}, or \cite[Section 1.5]{Fu}.
In this way, a polarized toric variety $(X, L)$ corresponds to an integral convex polytope $P$.
\end{remark}

\begin{definition}  An ample line bundle $L$ on a projective
variety $X$ is called {\it  normally generated} if the multiplication
map $\mbox{\textup{Sym}}^lH^0(X, L) \to H^0(X, L^{\otimes l})$ is
surjective for all $l\ge1$.  
\end{definition}

\begin{definition}
An integral convex polytope in $M_{\mathbb{R}}$ is called {\it normal} if
for the corresponding polarized toric variety $(X, L)$ the ample line bundle $L$ is
normally generated.
\end{definition}

\begin{remark}\label{rm}
If a polarized toric variety $(X, L)$ corresponds to
an integral convex polytope $P$ in $M_{\mathbb{R}}$, 
then the  normal generation of $L$ is equivalent to the normality of $P$.
From the result of Nakagawa \cite{N}, we see that it is enough to prove the equality
\begin{equation*}\label{sect1:eq0}
P\cap M + P\cap M = (2P)\cap M 
\end{equation*}
for the normality of $P$ in dimension three.
\end{remark}

In order to prove normality, we use the following easy lemma.

\begin{lemmama}\label{l1}
Let $P$ be a lattice polytope in $M_{\mathbb{R}}$.  If we have a cover $P=\cup_i Q_i$
by normal lattice polytopes $Q_i$, then $P$ is normal.
\end{lemmama}

\section{Proof of Theorem~\ref{int:t1}}\label{sect3}

In this section we prove  Theorem~\ref{int:t1}.  
Let $X$ be a nonsingular projective toric 3-fold and $A$ an ample line bundle on $X$
satisfying the condition that $\Gamma(A\otimes\mathcal{O}_X(K_X))=0$.
Let $Q$ be the integral convex polytope of dimension three corresponding to
the polarized toric variety $(X, A)$.  From  \cite[Theorem 3.6]{Od} we have
\begin{equation}\label{sect2:eq1}
\Gamma(X, A\otimes \mathcal{O}_X(K_X)) \cong 
\bigoplus_{m\in \mbox{{\scriptsize Int}}(Q)\cap M}
\mathbb{C}{\bold e}(m).
\end{equation}
Hence we see that $\Gamma(A\otimes\mathcal{O}_X(K_X))=0$ is equivalent to 
$\mbox{Int}(Q)\cap M=\emptyset$.

In Proposition 2.3 of \cite{Og2} Ogata classified nonsingular lattice polytopes $Q$ with 
$\mbox{Int}(Q)\cap M=\emptyset$.
From the classification, we know that $X$ is one of the following toric 3-folds.

\begin{itemize}
\item[(1)] A toric $\mathbb{P}^1$-bundle over a nonsingular toric surface $Y$.
\item[(2)] The projective 3-space $\mathbb{P}^3$.
\item[(3)] A $\mathbb{P}^2$-bundle over $\mathbb{P}^1$.
\item[(4)] The blowing up of $\mathbb{P}^3$ at torus invariant points (at most four points).
\item[(5)] The blowing up of a $\mathbb{P}^2$-bundle over $\mathbb{P}^1$ at one or two
torus invariant points such that two points do not lie on the same  fiber
 (at most two points).
\end{itemize}

When $X$ is the projective 3-space (the case (2)), 
set $P_1:=\mbox{Conv}\{(0,0,0), (1,0,0), (0,1,0), (0,0,1)\}$ the basic 3-simplex.  
Then $P_1$
 defines the polarized toric variety $(\mathbb{P}^3, \mathcal{O}(1))$.
 An ample line bundle on $\mathbb{P}^3$ is $\mathcal{O}(l)$ for some positive
 integer $l$.
 We know that $lP_1$ is normal for all $l\ge1$.
 
 In the cases (3) and (5), 
 $X$ has a surjective morphism onto the projective line.   From Theorem~\ref{keytm}, 
 all ample line bundles on $X$ are normally generated. 
  
 Before proceeding a proof, we introduce a useful lemma obtained in \cite[Lemma 2.3]{Og2}.
 This is also a corollary of the result of Fakhruddin\cite{Fa}.
 \begin{lemmama}\label{l2}
 Let $M'=\mathbb{Z}^2$.  Let $F, G$ be lattice polygons in $M'_{\mathbb{R}}$
 corresponding an ample and a nef line bundles on a nonsingular toric surface, respectively.
 Let $M:=M'\oplus \mathbb{Z}$.  Let $P=\mbox{\rm Conv}\{(F,0), (G,1)\}$ be the lattice polytope
 in $M_{\mathbb{R}}$ defined as the convex hull of $(F,0), (G,1)\subset M'_{\mathbb{R}}\times
 \mathbb{R}$.
 Then $P$ is normal. 
 \end{lemmama}
 
 Next, let $X$ be a toric $\mathbb{P}^1$-bundle over a nonsingular toric surface $Y$(the case (1)).
 Let $L$ be an ample line bundle in $X$ and $P$ the corresponding lattice polytope.
 $P$ has two special facets $F_0, F_1$ corresponding to ample line bundles on $Y$.
 By a suitable choice of coordinates $(x,y,z)$ in $M_{\mathbb{R}}$, we may assume that
 $P$ sits in the upper half space $(z\ge0)$, the faset $F_0$ is contained in the plane
 $(z=0)$ and $F_1$ is contained in $(z=n)$.  Set $G_i:=(z=i)\cap P$ for $i=0, \dots, n$.
 Since $P$ is nonsingular, $G_i$ are nonsingular polygons corresponding to 
 ample line bundles on $Y$.
 Set $P(G_i):=\mbox{conv}\{G_i, G_{i+1}\}$ for $i=0, \dots, n-1$.
 Then $P$ is the union of $P(G_i)$'s and each $P(G_i)$ is normal from Lemma~\ref{l2}.
 Thus $P$ is normal.
 
\begin{figure}[h]
 \begin{center}
 \setlength{\unitlength}{1mm}
  \begin{picture}(50,60)(15,5)
   \put(20,20){\vector(4,-1){30}}
   \put(20,20){\vector(0,1){30}}
   \put(20,20){\vector(2,1){33}}
   \put(2,15){\makebox(20,10){$(0,0,0)$}}
   \put(53,10){\makebox(10,10)[bl]{$x$}}
   \put(15,50){\makebox(10,10)[bl]{$z$}}
   \put(54,34){\makebox(10,10)[bl]{$y$}}
   \put(41.5,14){\line(2,1){11}}
   \put(42,14){\line(0,1){11}}
   \put(28,41){\line(5,-2){19}}
   \put(28,41){\line(2,-5){3}}
   \put(20,37){\line(4,-1){11}}
   \put(20,37){\line(2,1){8}}
   \put(47,33){\line(2,-5){5.5}}
   \put(42,25){\line(2,-1){10}}
   \put(42,25){\line(-5,4){12}}
   \put(35,5){\makebox(10,10)[r]{$(l_1,0,0)$}}
   \put(60,14){\makebox(10,10)[]{$(l_1,n-l_1,0)$}}
   \put(44,23){\makebox(10,10)[l]{$(l_1,0,n-l_1)$}}
   \put(43,37){\makebox(10,10)[br]{$(0,n,0)$}}
   \put(42,25){\circle*{1}}
   \put(47,33){\circle*{1}}
   \put(23,40){\makebox(10,10)[l]{$(0,n-m,m)$}}
   \put(21,23.5){\makebox(10,10)[tl]{$(n-m,0,m)$}}
   \put(2,32){\makebox(20,10){$(0,0,m)$}}
   \put(20,20){\circle*{1}}
   \put(42,14.5){\circle*{1}}
   \put(20,37){\circle*{1}}
   \put(52.5,19.5){\circle*{1}}
   \put(30.5,34.5){\circle*{1}}
   \put(28,41){\circle*{1}}
  \end{picture}
  \end{center}
  \caption{typical $P$ in case (4) of $l_2=l_3=0$}
 \label{fig0}
\end{figure}
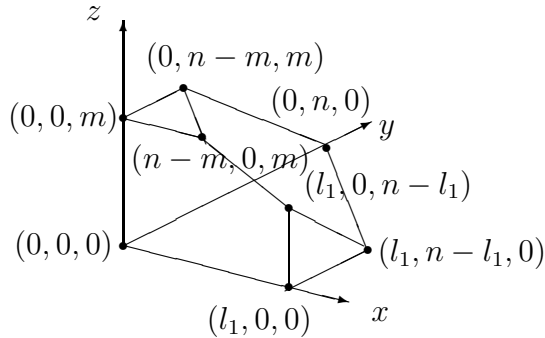

 In the last, let $X$ be at least one point blowup of $\mathbb{P}^3$(the case (4)).
 Let $L$ be an ample line bundle in $X$ and $P$ the corresponding lattice polytope.
 By a suitable choice of coordinates $(x,y,z)$ in $M_{\mathbb{R}}$, we may assume that
 $P$ sits in the region $(x\ge0, y\ge0,z\ge0)$ and three fasets $F_1, F_2, F_3$ are
 contained in the coordinate planes $(x=0),(y=0), (z=0)$, respectively.
 Since $P$ is obtained by several cuts from $nP_1$ for some positive integer $n$,
 $P$ has the facet $F_0$ contained in the plane $(x+y+z=n)$.
 By definition, we may assume that $P$ has the facet $F_4$ contained in the plane $(z=m)$
 for $1\le m\le n-1$.
 Moreover, $P$ may have facets $F_5, F_6, F_7$ contained in the planes $(x=l_1), (y=l_2),
 (x+y+z=l_3)$ for $m+1\le l_1, l_2\le n-1$ and $1\le l_3\le m-1$, respectively. 
 See Figure~\ref{fig0}.
 Set $G_i:=(z=i)\cap P$ for $i=0, \dots, m$.
 Since $P$ is nonsingular, each $G_i$ is a nonsingular polygon corresponding to an ample line bundle on a nonsingular toric surface $Y_i$.
 Set $P(G_i):=\mbox{conv}\{G_i, G_{i+1}\}$ for $i=0, \dots, m-1$.  Since $G_{i+1}$ is nef on $Y_i$,
 $P(G_i)$ is normal from Lemma~\ref{l2}.
 Since $P$ is the union of $P(G_i)$'s, $P$ is normal from Lemma~\ref{l1}.

 This completes the proof of Theorem~\ref{int:t1}.
 
 \section{Proof of Theorem~\ref{int:t2}}\label{sect4}
 
 In this section we prove  Theorem~\ref{int:t2}.  
Let $X$ be a nonsingular projective toric 3-fold and $A$ an ample line bundle on $X$
satisfying the condition that $\Gamma(A\otimes\mathcal{O}_X(K_X))\not=0$
and $A\otimes\mathcal{O}_X(K_X)$ is not big.
Let $Q$ be the integral convex polytope of dimension three corresponding to
the polarized toric variety $(X, A)$.   Then $\mbox{Int}(Q)\cap M\not=\emptyset$
and $ \dim\mbox{Conv}\{\mbox{Int}(Q)\cap M\}\le2$.

Set $Q^{\circ}:=\mbox{Conv}\{\mbox{Int}(Q)\cap M\}$.
We note that $Q^{\circ}$ also corresponds to a nef line bundle $B$ on $X$.
If $\dim Q^{\circ}=0$, then $B\cong \mathcal{O}_X$.

We separate the argument into three cases according to the dimension of $Q^{\circ}$.

The case I: $\dim Q^{\circ}=1$. Then $B$ defines the surjective morphism onto the projective 
line.  From Theorem~\ref{keytm},  all ample line bundles on $X$ are normally generated. 
 
The case II:  $\dim Q^{\circ}=2$.  We see that $Q$ has special two facets $F_0$ and $F_1$
parallel to $Q^{\circ}$.
Let $Y_0, Y_1$ be toric subvarieties of $X$ corresponding to $F_0, F_1$, respectively.
Set $H$ the plane in $M_{\mathbb{R}}$ containing $Q^{\circ}$ and set $G:=Q\cap H$.
Then $G$ is a lattice polytope of dimension two, and it corresponds to  nef and big line bundles
on the nonsingular toric surface $Y_0$ or $Y_1$.
If $G$ is ample on both $Y_0$ and $Y_1$, then $X$ has a structure of $\mathbb{P}^1$-bundle
over $Y_0$.  We treated this case in the previous section as the case (1).
If $G$ is not ample on $Y_0$, then we see that $G$ corresponds to an ample line bundle on
a nonsingular toric surface $Y$ obtained by blowing down $(-1)$-curves on $Y_0$
from the argument in \cite{Og2}.
From the theory of minimal models of rational surfaces, $Y$ has a surjective morphism onto 
the projective line unless $Y\cong \mathbb{P}^2$.
In this case $X$ has also a surjective morphism onto the projective line.
We can apply Theorem~\ref{keytm}.

We assume $Y\cong \mathbb{P}^2$.
Then $X$ is blowup at most three invariant points of a $\mathbb{P}^1$-bundle over $\mathbb{P}^2$.
We note that any two blowup points are not contained in the same fiber.
Let $Y_0, Y_1$ be toric subvarieties of $X$ blowup  $r, s$ points ($r\ge s$) of $\mathbb{P}^2$, respectively.
Let $L$ be an ample line bundle on $X$ and $P$ the lattice polytope corresponding to $L$.
Set $F_i$ the facet of $P$ corresponding to $L|_{Y_i}$ for $i=0,1$.
By a suitable choice of coordinates $(x,y,z)$ in $M_{\mathbb{R}}$, we may assume that
 $P$ sits in  the region $(x\ge0, y\ge0,z\ge0)$, the facet $F_0$ is contained in the plane
 $(z=0)$, $F_1$ is contained in $(z=n)$ and $P$ has two facets contained in the coordinate planes $(x=0), (y=0)$.  
 Set $G_i:=(z=i)\cap P$ for $i=0, \dots, n$.
Set $P(G_i):=\mbox{conv}\{G_i, G_{i+1}\}$ for $i=0, \dots, n-1$.

If $s=0$, then $F_1=G_n$ is a triangle.  Set $Z_i$ the toric surface corresponding to the
lattice polygon $G_i$ for $i=0, \dots, n$.  Then $G_{i+1}$ is nef on $Z_i$.
Thus each $P(G_i)$ is normal from Lemma~\ref{l2}.

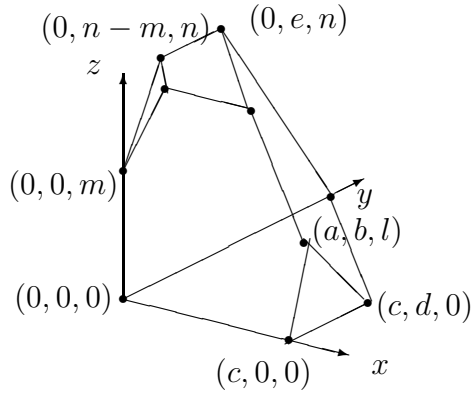
\begin{figure}[h]
 \begin{center}
 \setlength{\unitlength}{1mm}
  \begin{picture}(50,60)(15,5)
   \put(20,20){\vector(4,-1){30}}
   \put(20,20){\vector(0,1){30}}
   \put(20,20){\vector(2,1){32}}
   \put(2,15){\makebox(20,10){$(0,0,0)$}}
   \put(53,10){\makebox(10,10)[bl]{$x$}}
   \put(15,50){\makebox(10,10)[bl]{$z$}}
   \put(51,32){\makebox(10,10)[bl]{$y$}}
   \put(41.5,14){\line(2,1){11}}
   \put(42,14){\line(1,5){2.8}}
   \put(33,56){\line(2,-3){15}}
   \put(33,56){\line(1,-3){3.5}}
   \put(25,52){\line(1,-5){0.9}}
   \put(26,48){\line(4,-1){11}}
   \put(25,52){\line(2,1){8}}
   \put(20,37){\line(1,2){5.5}}
   \put(20,37){\line(1,3){5}}
   \put(47.5,33.5){\line(2,-5){5.5}}
   \put(37,45){\line(2,-5){7}}
   \put(44,28){\line(1,-1){8}}
   \put(35,5){\makebox(10,10)[r]{$(c,0,0)$}}
   \put(55,14){\makebox(10,10)[]{$(c,d,0)$}}
   \put(45,24){\makebox(10,10)[l]{$(a,b,l)$}}
   \put(40,55){\makebox(10,10)[br]{$(0,e,n)$}}
   \put(44,27.5){\circle*{1}}
   \put(47.6,33.6){\circle*{1}}
   \put(9,51){\makebox(10,10)[l]{$(0,n-m,n)$}}
   \put(2,30){\makebox(20,10){$(0,0,m)$}}
   \put(20,20){\circle*{1}}
   \put(42,14.5){\circle*{1}}
   \put(20,37){\circle*{1}}
   \put(52.5,19.5){\circle*{1}}
   \put(33,56){\circle*{1}}
   \put(25,52){\circle*{1}}
   \put(25.5,48){\circle*{1}}
   \put(37,45){\circle*{1}}
  \end{picture} \end{center}
\caption{typical $P$ with $s=r=1$}
 \label{fig1}
\end{figure}

We assume $s=1$.  By assumption we have $1\le r\le 2$.
Set $1\le m<n$.  We assume that $P$ has a facet $F_2$ contained in the plane $(z=x+y+m)$.
See Figure~\ref{fig1}.
Since $1\le r\le2$, $P$ has at least one and may have two facets $F_3, F_4$ 
corresponding the blowup 
$Y_0\to \mathbb{P}^2$.
Set $l_j$ the largest $z$-coordinate of vertices of $F_j$ for $j=3,4$.
Set $l_3\le l_4$.

For $0\le i\le m-1$ or $m+1\le i<n$, the line bundle corresponding $G_{i+1}$ on $Z_i$
is nef, hence, $P(G_i)$ is normal.
Even if $i=m$, $G_{m}$ is nef on $Z_{m+1}$ unless $l_4\le m$ or $l_3\ge m+2$.

We have to consider the cases $l_3<l_4=m+1$, $l_3=l_4=m+1$ and $l_3=m+1<l_4$.
We may imagine the shape of $P(G_m)$ by setting $m=0, n=l=1$ in Figure~\ref{fig1}.
We note that $d=1, c\ge1, a, e\ge2$ in the Figure.
We decompose $P(G_m)$ as the union of $R_+:=P(G_m)\cap(z\ge x+y+m)$ and
$R_-:=P(G_m)\cap(z\le x+y+m)$.
Set $G_i^{\pm}:=G_i\cap R_{\pm}$.  From Lemma~\ref{l2}, $R_{\pm}$ are normal because
$G_m^+$ is nef with respect to $G_{m+1}^+$ and $G_{m+1}^-$ is nef with respect to
$G_m^-$.
Thus $P$ is normal from Lemma~\ref{l1}.
  
The case III: $\dim Q^{\circ}=0$.
We recall the situation when $A\otimes\mathcal{O}_X(K_X)$ is not nef in \cite{Og2} for an ample line bundle $A$ on
a nonsingular toric 3-fold $X$.
When $\Gamma(X, A\otimes\mathcal{O}_X(K_X))\not=0$, 
if $A\otimes\mathcal{O}_X(K_X)$, then there exists an ample line bundle $\Bar{A}$ on a nonsingular toric 3-fold $Y$ 
 and a surjective morphism $\pi: X \to Y$ which is blow up at several points such that
 $\Bar{A}\otimes\mathcal{O}_Y(K_Y)$ is nef and $\Gamma(X, A\otimes\mathcal{O}_X(K_X))\cong 
 \Gamma(Y, \Bar{A}\otimes\mathcal{O}_Y(K_Y))$.
 Here if we set $\sum_i E_i$ the exceptional divisor of $\pi$, then 
 $A=\pi^*\Bar{A}\otimes\mathcal{O}_X(-\sum_i E_i)$, $K_X=\pi^*K_Y\otimes\mathcal{O}_X(
 \sum_i2E_i)$ and $A\otimes\mathcal{O}_X(K_X-\sum_i E_i)=\pi^*(\Bar{A}\otimes\mathcal{O}_Y(K_Y))$.
 Since in our case $B=\pi^*(\Bar{A}\otimes\mathcal{O}_Y(K_Y))\cong \mathcal{O}_X$, we see that $Y$ is a toric
 Fano 3-fold.

Batyrev\cite{B} and Watanabe-Watanabe\cite{WW} classified nonsingular toric Fano 3-folds.
By the classification we see that there exists 18 toric Fano's up to isomorphism.
If $Y$ has a surjective morphism onto the projective line, then so does $X$, hence,
all ample line bundles on $X$ are normally generated from Theorem~\ref{keytm}.
According to the classification of Oda's book\cite{Od}, except the class admitting surjective
morphism onto the projective line, we have four classes, that is,
(1) the $\mathbb{P}^3$, (3) and (4) 
two $\mathbb{P}^1$-bundle
over $\mathbb{P}^2$ and (11).

 Set $\Bar{Q}$ the lattice polytope corresponding $\mathcal{O}_Y(-K_Y)$ of a toric Fano 3-fold $Y$.
For the class (1), $\Bar{Q}=4P_1$.  Thus we see that $X$ is blowup at most four points of
$\mathbb{P}^3$.  This is the case (1) in Section~\ref{sect3}. 
 
For the class (3), $\Bar{Q}=(4P_1)\cap(0\le x\le2)$.  Thus we see that $X$ is blowup at most
 three points of $Y$ whose any two points are not contained in the same fiber.
 This is the case II in this section.
 
 For the class (4), $\Bar{Q}=\mbox{Conv}\{(0,0,0), (5,0,0), (0,5,0), (0,0,2), (1,0,2), (0,1,2)\}$.
 Thus we see that $X$ is blowup at most
 three points of $Y$ whose any two points are not contained in the same fiber.
 This is also the case II in this section.

For the class (11), $\Bar{Q}=\mbox{Conv}\{(0,0,0), (4,0,0), (0,4,0), (0,0,2), (1,0,2), (0,1,2),
(3,0,1), (0,3,1)\}$.
Thus we see that $X$ is blowup at most one point of $Y$.
Let $P$ be a integral convex polytope corresponding to an ample line bundle on this $X$. 
$P$ has two parallel facets $F_0, F_1$.  By cut at lattice hyperplanes parallel to $F_0$,
we can decompose $P$ into a union of normal lattice polytopes as in the case II in this section.
Thus we see that $P$ is normal from Lemmas~\ref{l2} and \ref{l1}.

This completes the proof of Theorem~\ref{int:t2}.

\bigskip
\address{
Mathematical Institute \\ 
Tohoku University \\
Sendai 980-8578 \\
Japan
}
{ogata@math.tohoku.ac.jp}

\end{document}